\magnification=\magstep1
\hsize=32pc
\vsize=43pc
\parindent 0cm

\input epsf.tex


\font\bigbf=cmbx10 scaled\magstep3

\def\ni{\noindent}
\def\eop{{ \vrule height7pt width7pt depth0pt}\par\bigskip}
\def\NN{{ I\!\!N}}

\def\QQ{{\rlap {\raise 0.4ex \hbox{$\scriptscriptstyle |$}}
\hskip -0.2em Q}}
\def\RR{{ I\!\!R}}

\def\ZZ{{ Z\!\!\! Z}}
\def\CC{{ \rlap {\raise 0.4ex \hbox{$\scriptscriptstyle |$}} \hskip -0.2em C}}

\def\s{{\bf s}}
\def\t{{\bf t}}
\def\c{{\bar c}}


\centerline{\bigbf Limit laws of entrance times}
\bigskip
\centerline{\bigbf for low complexity Cantor minimal systems} 
\vskip 1.5cm
\centerline{\bf Fabien Durand} 
\bigskip
\ni Facult\'{e} de Math\'{e}matiques et
d'Informatique et Laboratoire Ami\'enois
de Math\'ematiques Fondamentales  et
Appliqu\'ees, CNRS-ESA 6119, Universit\'{e} de Picardie
Jules Verne, 33 rue Saint Leu, 80000 Amiens, France,
and Centro de Modelamiento Matem\'atico, \hfill\break UMR 2071 UCHILE-CNRS,
e-mail:fabien.durand@u-picardie.fr, fdurand@dim.uchile.cl
\bigskip
\centerline{\bf Alejandro Maass} 
\bigskip
\ni Departamento de Ingenier\'{\i}a Ma\-te\-m\'a\-ti\-ca and
Centro de Modelamiento Matem\'atico, \hfill\break UMR 2071 UCHILE-CNRS,
Universidad de Chile, Facultad de Ciencias F\'{\i}sicas y
Mate\-m\'a\-ti\-cas, Casilla 170-3, Correo 3, Santiago, Chile.
e-mail: amaass@dim.uchile.cl
\vskip 1.5cm
\centerline{\bf Abstract}
\medskip
\ni {\it 
This paper is devoted to the study of limit laws 
of entrance times to cylinder sets 
for Cantor minimal systems of zero entropy 
using their representation 
by means of ordered Bratteli diagrams.
We study in detail 
substitution subshifts and  we prove  
these  limit laws are piecewise linear functions. 
The same kind of results is obtained for  classical
low complexity systems given by  non stationary ordered
Bratteli diagrams.  }
\bigskip
\bigskip
\ni{\bf 1. Introduction.}
\bigskip
\bigskip
\ni{\it 1.1. Preliminaries and motivations.}
\bigskip
A topological dynamical system, or just dynamical system, is a compact
Hausdorff space $ X $ together with a homeomorphism $ T:X\rightarrow X.$ 
We denote it by  $ \left( X,T\right)  $. If $X$ is a Cantor set we
say that $\left( X,T\right)$  is a Cantor system. That is, $ X $ has a countable basis of closed
and open sets (clopen sets) and it has no isolated points. 
A dynamical system is minimal
if all orbits $\{T^n(x):n\in \ZZ\}$ are dense in $X$, or equivalently the only non trivial closed
$T$-invariant set is $X.$ 
\medskip
Let $(X,T)$ be a Cantor minimal system and fix a $T$-invariant probability measure
$\mu$. Let $I\subseteq X$ be a clopen set.
For each $x \in X$ the entrance time to $I$ and the $k$-th return time
to $I$ for $k\ge 2$ are defined respectively by
$$N^{(1)}_I(x)=\inf \{ n>0: T^n(x) \in I \} \hbox{ and }
N^{(k)}_I(x)=\inf \{ n> N_I^{(k-1)}(x): T^n(x) \in I \}.$$
Since the system is minimal these quantities are finite.
The corresponding distributions are
$$F_I^{(1)}(t)=\mu\{x \in X : \mu(I) \cdot N^{(1)}_I(x) \le t \},$$
and, for $k>1$,
$$F_I^{(k)}(t)=\mu\{x \in X:
\mu(I) \cdot (N^{(k)}_I(x)- N^{(k-1)}_I(x)) \le t \} .$$
\medskip
Consider  the following problem:
fix a point $x^* \in X$
and let $\mu$ be a $T$-invariant probability 
measure of $(X,T)$. Let $(I_n:n\in \NN)$ be  a sequence of clopen sets
of $X$ such that $x^* \in I_n$,
$I_{n+1}\subseteq I_n$ for all $n \in \NN$, and
$\cap_{n \in \NN} I_n=\{x^*\}$.
For each $x \in X$ and every $k\ge 2$ define
$$N^{(1)}_n(x)=N^{(1)}_{I_n}(x)  \hbox{ and } N^{(k)}_n(x)=N^{(k)}_{I_n}(x).$$
We will study the limits of the
sequences of distributions $(F^{(1)}_{I_n})_{n\in \NN}$ and
$(F^{(k)}_{I_n})_{n\in \NN}$ for $k\ge 2$ and $(I_n:n\in \NN)$ a sequence of 
cylinder sets. 
For simplicity we will write $( F^{(1)}_n )_{n\in \NN}$ and $( F^{(k)}_n )_{n\in \NN}$ respectively.
These limits, when they exist,  will be called limit laws of entrance times.
\medskip
The existence and characterization of limit laws 
for particular (and natural) families of sequences $(I_n:n\in \NN)$
is a problem that has been addressed  in several 
papers in the last ten years. 
Most of them has focused on systems of positive entropy with
strong conditions of mixing (see [CC1,CC2,CG,H,HSV,P]).
In all of these cases the limit laws are exponential.
The unique non exponential limit laws we know appear in the 
study of homeomorphisms of the circle [CdF]. 
Under some mild  conditions on the 
continued  fraction expansion of the  rotation numbers 
the authors found piecewise linear limit laws. 
In this work $(I_n:n\in \NN)$ is a sequence of intervals which 
end points are given by the partial quotients of the continued 
fraction expansion of the angle.
They 
proved under the same assumptions the convergence in law  
of the associated point process. 
\medskip
The present paper is motivated by our reading  of [CdF]. In this work the main
arguments concern irrational rotations of the interval. 
These systems are measure--theoretically conjugate to Sturmian subshifts introduced in
[HM].  In the symbolic context they correspond to non trivial subshifts with the 
lowest complexity. Also, we know that whenever the rotation number is quadratic 
the associated Sturmian subshift is a substitutive subshift [DDM].
In [C] the author addressed the question whether 
analog results as those in [CdF] could appear in the context of
substitutive subshifts. Moreover, the author  expected that weak mixing  
would be necessary.
\medskip
In the present work we address the same questions described before  
in the  
framework of minimal Cantor systems, represented 
by means of Bratteli--Vershik systems, for  
sequences $(I_n:n\in \NN)$ made of cylinder sets. 
In particular, we provide answers 
for substitution subshifts, odometers and Sturmian sequences.
In all these cases we get under some mild assumptions 
piecewise linear limit laws (Theorem 2.4, Examples 3,4,5).
The main tool developed here is a counting procedure
over the ordered Bratteli diagrams used to represent those systems. 
The representation of Cantor minimal systems 
by means of  ordered Bratteli diagrams has been 
introduced in [HPS] and it has been used to solve the problem of orbit equivalence (we present them
below).
Nowadays, there exist characterizations of these diagrams
for large classes of subshifts; in particular, substitution subshifts [DHS],
Sturmian subshifts [DDM], linearly recurrent subshifts
[D]  and Toeplitz subshifts [GJ]. 
The nice structure of such diagrams allows to reduce most of the
problems to a matrix analysis.
Finally, in section 4 we study the point process associated to 
entrance times of substitution subshifts. We point out that we never assume 
any mixing condition. 
\medskip
By the time we were submitting this article 
Y. Lacroix [L] has obtained the following  general result:  
given an aperiodic ergodic system $(X,{\cal B},\mu,T)$ and
a distribution function $G:\RR \to [0,1]$ 
there exists a sequence $(I_n:n\in \NN)$ 
such that 
$\mu\{x \in I_n: \mu(I_n) N_{I_n}^{(1)}(x) \le t \} / \mu(I_n) \to G(t)$ 
as $n\to \infty$. This result used an  abstract construction 
based on Rokhlin towers. On the other hand the author provides
an explicit example of Toeplitz subshift where the sequence $(I_n:n\in \NN)$
consists of cylinder sets.
\bigskip
\ni{\it 1.2. Subshifts and Complexity. }
\bigskip
A particular class of Cantor systems is the class of 
subshifts. These systems are defined as follows. Take a finite set or
alphabet $A$. The set $ A^{\ZZ }$ 
consists of infinite sequences $ \left( x_{i}\right) _{i\in \ZZ }$ 
with coordinates $ x_{i}\in A. $ With the product topology $ A^{\ZZ } $ 
is a compact Hausdorff Cantor space. We define the shift transformation
$ \sigma :A^{\ZZ }\rightarrow A^{\ZZ } $
by $ \left( \sigma \left( x\right) \right)_{i}=x_{i+1}$ 
for any $ x\in A^{\ZZ } $, $ i\in \ZZ  $. 
The pair $ \left( A^{\ZZ },\sigma \right)$ 
is called a full shift. A subshift is a pair $(X,\sigma )$ where
$X$ is any $\sigma$--invariant closed subset of $A^{\ZZ }$.
A classical procedure to construct subshifts is by considering the closure of
the orbit under the shift of a single sequence 
$ x\in A^{\ZZ } $ , 
$\Omega (x)= \overline{\left\{ \left. \sigma ^{i}\left( x\right) 
\ \right| \ i\in \ZZ \ \right\} } $.
\medskip
Let $  \left( x_{i}\right) _{i\in \NN } $  be an element of $  A^{\NN } $.
Another classical procedure is to consider the set 
$ \Omega \left( x\right) $ 
of infinite sequences $ \left( y_{i}\right) _{i\in \ZZ } $ such that for
all $ i\leq j $ there exists 
$k\geq 0 $ such that 
$y_{i} y_{i+1}\cdots y_{j}=x_{k}x_{k+1}\cdots x_{k+j-i}$.
In both cases we say that 
$ \left( \Omega \left( x\right) ,\sigma \right) $ 
is the subshift generated by $ x $.
\medskip
A classical measure of complexity of a zero entropy subshift  
$(X,T)$  
is the so called symbolic complexity. It is the integer function
$p_X:\NN \to \NN$ where $p_X(n)$ is the number of all different 
words of length $n$ appearing in sequences of $X$. We say that 
the complexity is sub-linear if there exists a positive constant $a$ 
such that $p_X(n)\le a n$.
\bigskip
\ni{\it 1.3. Bratteli--Vershik representations. }
\bigskip
A Bratteli diagram is an infinite graph $\left( V,E\right)$ which
consists of a vertex set $V$ and an edge set $E$, both of which are
divided into levels $V=V_{0}\cup V_{1}\cup \cdots  $, 
$ E=E_{1}\cup E_{2}\cup \cdots $ 
and all levels are pairwise disjoint. The set $ V_{0} $ is a 
singleton $ \{v_{0}\} $,
and for $k\geq 1$, $E_{k}$  is the set of edges joining vertices in
$V_{k-1}$ to vertices in $V_{k}$. 
It is also required that every vertex
in $V_{k}$ is the ``end-point'' of some edge in $ E_{k} $ for $ k\geq 1 $,
and the ``initial-point'' of some edge in $E_{k+1} $ for $ k\geq 0 $.
By level $k$ we will mean the subgraph consisting of the
vertices in $V_{k}\cup V_{k+1}$ and the edges $E_{k+1}$ between these
vertices. 
We describe the edge set $E_{k}$ using a $V_{k-1}\times V_{k}$ incidence
matrix, $M^{(k)}$, for which its $(i,j)$--entry is the number of edges in $E_{k}$
joining vertex $i\in V_{k-1}$ with vertex $j\in V_{k}$.
For every $e\in E_k$,  $\s (e)\in V_{k-1}$ and $\t (e) \in V_k$ are the starting and terminal vertices 
of $e$ respectively.
\medskip
An ordered Bratteli diagram $ B=\left( V,E,\preceq \right)  $ is a
Bratteli diagram $ \left( V,E\right)  $ together with a partial ordering
$ \preceq  $ on $ E $. Edges $ e $ and $ e' $ are comparable if
and only if they have the same end-point. We call $succ(e)$ the successor 
of $e$ with respect to this partial order when $e$ is not a maximal edge.
\medskip
\centerline{}
\epsfxsize=12truecm
\medskip
\centerline{\epsfbox{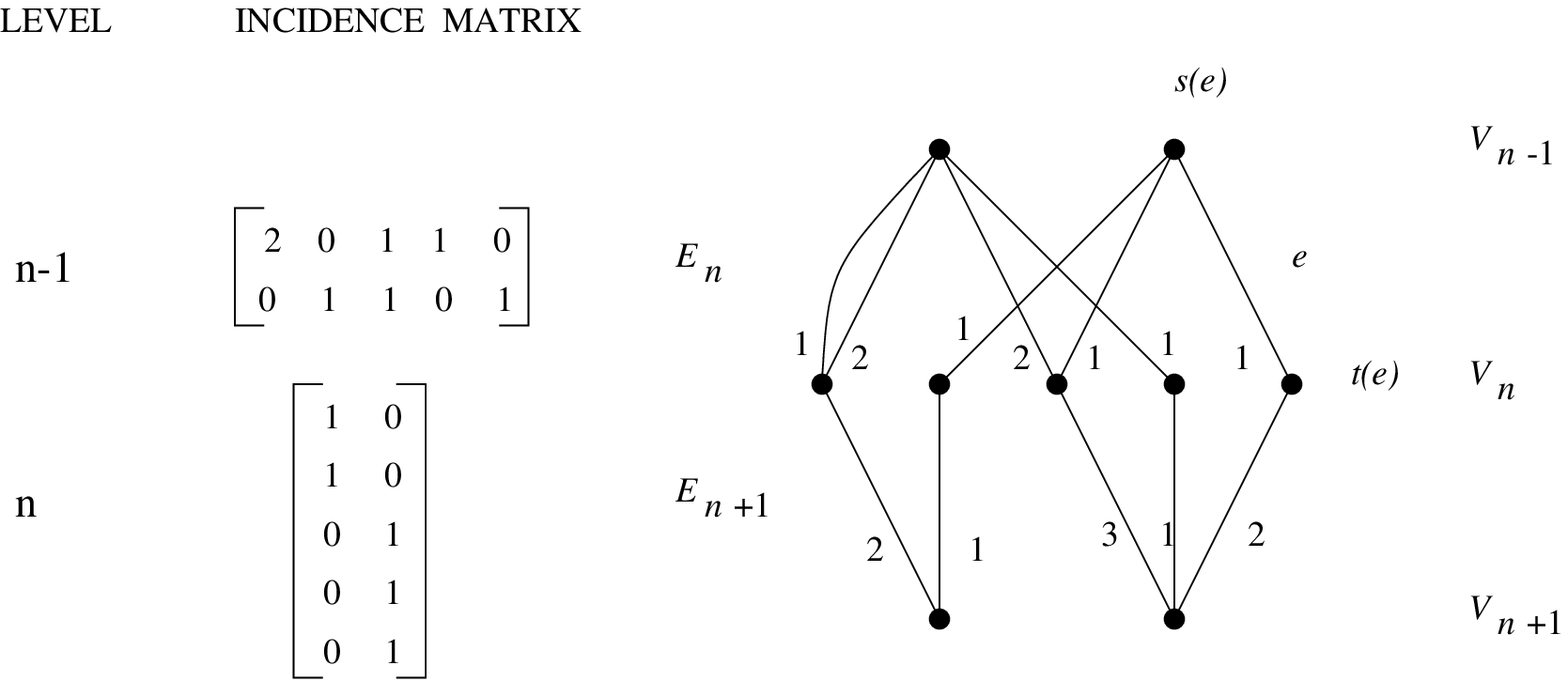}}
\medskip
\centerline{\bf Figure 1}
\medskip
Let $k<l$ in $\NN \setminus \left\{ 0\right\}$ and let $E(k,l)$
be the set of all paths of length $l-k$ 
in the graph joining vertices of $V_{k-1}$ with
vertices of $V_{l}$. The partial ordering of $E$ induces another in
$E(k,l)$ given by $\left( e_{k},\ldots ,e_{l}\right) 
\prec \left( f_{k},\ldots ,f_{l}\right) $ 
if and only if there is $ k\leq i\leq l $ such that $ e_{j}=f_{j} $ for
$ i<j\leq l $ and $ e_{i}\prec f_{i}$.
\medskip
Given a strictly increasing sequence of integers 
$ \left( m_{n}\right) _{n\geq 0}$ 
with $m_{0}=0$ we define the contraction 
of $ B=\left( V,E,\preceq \right) $ 
(with respect to $ \left( m_{n}\right) _{n\geq 0} $) as 
$$\left( \left( V_{m_{n}}\right) _{n\geq 0},
\left( E(m_{n}+1,m_{n+1})\right) _{n\geq 0},\preceq \right),  $$
where $ \preceq  $ 
is the order induced in each set of edges 
$ E(m_{n}+1,m_{n+1}) $. The inverse operation of contracting is microscoping (see [GPS]).
\medskip
We say that an ordered Bratteli diagram is stationary if 
for any $k\geq 1$ the incidence matrix and order are the same 
(after labeling the vertices appropriately).
\medskip
Given an ordered Bratteli diagram 
$B=\left( V,E,\preceq \right)$ we define
$X_{B}$ as the set of infinite paths 
$\left( e_{1},e_{2},\cdots \right)$ 
starting in $v_{0}$ such that for all 
$i\geq 1$ the end-point of $e_{i}\in E_{i}$
is the initial-point of 
$e_{i+1}\in E_{i+1}$. 
We topologize $X_{B}$
by postulating a basis of open sets, namely the family of 
cylinder sets
$$
[e_{1},e_{2},\ldots ,e_{k}] =\left\{  \left( 
f_{1},f_{2},\ldots \right) \in X_{B}
\ : \ f_{i}=e_{i},\hbox{ for } 1\leq i\leq k\ \right\} .$$
Each $[e_{1},e_{2},\ldots ,e_{k}]$ is also closed, as is
easily seen, and so we observe that 
$X_{B}$ is a compact, totally disconnected
metrizable space.
\medskip
When there is a unique $x=\left( x_{1},x_{2},\ldots \right) \in X_{B}$
such that $x_{i}$ is maximal for 
any $i\geq 1$ and a unique $y=\left( y_{1},y_{2},\ldots \right) \in X_{B}$
such that $y_{i}$ is minimal 
for any $i\geq 1$, we say that $B=\left( V,E,\preceq \right)$
is a properly ordered Bratteli diagram. Call these particular points
$x_{\hbox{max}}$ and $x_{\hbox{min}}$ respectively. In this case
we can define a dynamic $V_{B}$ over $X_{B}$ called Vershik
map. The map $V_{B}$ 
is defined as follows: let $x=\left( e_{1},e_{2},\ldots \right) 
\in X_{B}\setminus \left\{ x_{\hbox{max}}\right\} $ 
and let $k\geq 1$ be the smallest integer so that $e_{k}$ is not a
maximal edge. Let $f_{k}$ 
be the successor of $e_{k}$ and $\left( f_{1},\ldots ,f_{k-1}\right)$ 
be the unique minimal path in $E_{1,k-1}$  connecting $v_{0}$  with
the initial point of $f_{k}$. 
We set 
$V_{B}\left( x\right) =
\left( f_{1},\ldots ,f_{k-1},f_{k},e_{k+1},\ldots \right) $ 
and 
$V_{B}\left( x_{\hbox{max}}\right) =x_{\hbox{min}}$. The dynamical
system $\left( X_{B},V_{B}\right)$ is called Bratteli-Vershik system
generated by $B=\left( V,E,\preceq \right)$. 
The dynamical system induced
by any contraction of $B$ 
is topologically conjugate to $\left( X_{B},V_{B}\right)$.
In [HPS] it is proved that any minimal Cantor system 
$\left( X,T\right)$
is topologically conjugate to a Bratteli-Vershik system 
$\left( X_{B},V_{B}\right)$.
We say that $\left( X_{B},V_{B}\right)$ is a Bratteli-Vershik representation
of $\left( X,T\right)$.
\bigskip
\bigskip
\ni{\bf 2. Limit laws for stationary Bratteli--Vershik systems.}
\bigskip
\bigskip
Let us begin with some additional definitions and background.
A substitution is a map $\tau :A\rightarrow A^{+}$, where $A^{+}$
is the set of finite sequences with values in
$ A $. We associate to $\tau$
a $A\times A$ square matrix $M_{\tau }=\left(m_{a,b}\right)_{a,b\in A}$
such that $m_{a,b}$ is the number of times that the letter $ a $
appears in $ \tau \left( b\right)  $.
We say that $ \tau  $ is primitive
if $M_{\tau }$ is primitive, i.e. if some power of $M_{\tau }$ has
strictly positive entries only. A substitution $\tau$ can be naturally
extended by concatenation to $A^{+}$, $A^{\NN }$ and $A^{\ZZ }$.
We say that a subshift of
$A^{\ZZ }$ is generated by the substitution $\tau$
if it is the orbit closure of a fixed point for $\tau$ in $A^{\NN }$. It is well known that
primitivity of $\tau$ implies that this subshift  is minimal and uniquely ergodic
(see [Q] for more details).
\medskip
Let $(p_k:k\in \NN)$ be a sequence of positive integers.
The inverse limit of the sequence of groups $(\ZZ / p_1\cdots p_k \ZZ: k \in \NN)$
endowed with the addition of 1 is called odometer with base $(p_k:k\in \NN)$.
These systems are minimal and uniquely ergodic.
We say it is of constant base if 
the sequence $(p_k:k\in \NN)$ is ultimately constant.
\medskip
In [DHS] (see also [F]) it is proved that the family of stationary  
Bratteli--Vershik systems is up to topological conjugacy 
the disjoint union of the family of  
substitution minimal subshifts and the family of odometers with constant base. 
\medskip
Let $(X_B,V_B)$ be the  minimal Cantor system given by 
the  stationary ordered Bratteli diagram 
$\displaystyle {B}=(\cup_{i\ge 0} V_i, \cup_{i\ge 1} E_i, \preceq)$ where       
$V_i=\{ v(i,1),...,v(i,m) \}$, for  $i\ge 1$, and $V_0=\{v_0\}$. 
Moreover, by an appropriate labeling of the vertices 
the incidence matrices $(M^{(i)}:i\ge 1)$ are all equal to a matrix $M$.
In the sequel we identify each $V_i$ to $\{1,..,m\}$ following the labeling 
of vertices chosen to define $M$. In this setting the order of edges 
is the same for any level greater than one.
This representation is not unique and 
in this paper we will consider one that is appropriate for our purpose.  In the sequel 
we fix one which satisfies:
\medskip
(H1) the incidence matrix, $M$,  of $B$ has strictly positive coefficients;
\smallskip
(H2) for every vertex  $i\in V_1$ there is a unique edge from $v_0$ to $i$;
\smallskip
(H3) $\forall i \in \{1,...,m\}$, $\forall n \ge 1$,
$e=min\{f \in E_n:\t(f)=i\} \Rightarrow \s(e)=1$; 
\medskip
Let us notice that this representation can always be obtained 
contracting and microscoping levels if necessary. 
We recall that for  all $n \ge 2$,  
$M_{i,j}=|\{e \in E_n: \s (e)=i, \t (e)=j \}|$ and we also remark that
$\sum_{i=1}^m M^{n-1}_{i,j}$ is
the number of paths of length $n$ joining  $v_0$ with $j \in V_n$.
\medskip
Let $\lambda$ be the 
maximal  eigenvalue of $M$. We denote by
$r=(r(i):i\in \{1,...,m\})^T$ and
$l=(l(i):i\in \{1,...,m\})$ the
corresponding strictly positive right and left eigenvectors respectively,
such that $\sum_{i=1}^mr(i)=1$, $\sum_{i=1}^ml(i)\cdot r(i)=1$.
For every $e_1...e_n\in E(1,n)$,  we have that the unique ergodic measure is
defined by
$\mu([e_1...e_n])={r(\t(e_n)) \over \lambda^{n-1}}$ (for more details 
on the construction of measures for Bratteli-Vershik systems 
you can in particular see [BJKR]).
\medskip
\ni{\bf Example 1:} Consider the system given by the Bratteli
diagram in Figure 2. The order is written  over the edges and 
the  incidence associated matrix is $M=\left [ \matrix{ 1 & 1 \cr 2 & 3 \cr } 
\right ]$. \eop
\medskip
\epsfxsize=4truecm
\medskip
\centerline{\epsfbox{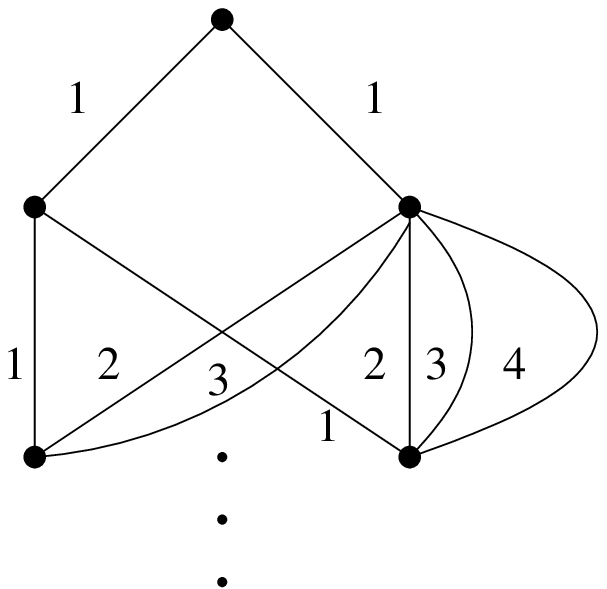}}
\medskip
\centerline{\bf Figure 2}
\bigskip

Let us now present the main problem of the section.
Let $x^*=(x_1^*,x_2^*,...) \in X_B$ and consider
the cylinder sets induced by $x^*$, that is 
$I_n=I_n(x^*)=[x^*_1,...,x^*_n]$.
We will study the limit laws of entrance times 
for this family of cylinder sets.
\medskip
Since the diagram $B$  
is stationary there is $i^* \in \{1,...,m\}$ such that 
$\t(x^*_{n})=i^*$ infinitely often.
Let ${\cal N}=(n_i:i \in \NN)$ be a subsequence such that
$\t(x^*_{n_i})=i^*$. 
In order to compute the limit laws with respect to these subsequences, we need to know
$\mu\{x\in X_B: N_n^{(1)}(x)=j \}$ and 
$\mu\{x\in X_B:(N^{(k)}_n(x)- N^{(k-1)}_n(x))=j \}$ for $k\ge 2$, $j  \in 
\NN \setminus\{0\}$.
\bigskip
\ni{\bf Lemma 2.1} {\it Let $n \in {\cal N}$.

\ni (i) If  $ef\in E(n+1,n+2)$ with $\s(e)=i^*$, 
then for any $z,y \in [I_nef]$ we have  $N_n^{(1)}(z)=N_n^{(1)}(y)$.
\smallskip
\ni (ii) Let $k\ge 0$. If  $e_1...e_{k+1} \in E(n+1,n+k+1)$ with $\s(e_1)=i^*$, then there is 
$e_1^{(k)}e_2^{(k)} \in E(n+1,n+2)$ with $\s(e_1^{(k)})=i^*$  such that for any 
$z,y \in [I_ne_1...e_{k+1}]$, $N_n^{(k)}(z)=N_n^{(k)}(y)$ and 
$V_B^{N_n^{(k-1)}(z)}(z) \in [I_n e_1^{(k)} e_2^{(k)}],
V_B^{N_n^{(k-1)}(y)}(y)\in [I_ne_1^{(k)}e_2^{(k)}]$.}
\smallskip
\ni{\bf Proof:}

\ni (i) We observe that  
any point of the system belonging to the cylinder set generated by  
the minimal path from $v_0$ to any $i \in V_n$, moves  under the action 
of $V_B$, from this cylinder set to the one corresponding 
to the maximal path from $v_0$ to $i\in V_n$, passing 
successively, and respecting the order, by all the paths from 
$v_0$ to $i \in V_n$. Consequently the return times to $I_n$, 
that is $N^{(1)}_n(x)$ for $x \in I_n$, 
are the same as those computed for the minimal path connecting $v_0$ with 
$i^* \in V_n$. Let us call $J_n$ the minimal path 
joining $v_0$ with $i^* \in V_n$.
The dynamics of any point of $[J_nef]$ is the 
following (you can see Figure 3): 
(1) they move  from $[J_nef]$ to the maximal path from $v_0$
to $\t(e)$;
(2) they move  from this maximal path to the minimal path from $v_0$ to some 
$i' \in V_{n+1}$, where $i'=\s(succ(f))$ if $f$ is not maximal and $i'=1$ if it 
is maximal (condition (H3));
(3) finally, 
since there is an edge $e'$ from $i^*$ to $i'$ (condition (H1)),  
they move  from this minimal path to the maximal path connecting $v_0$
and $i'$ passing through the cylinder set $[J_ne']$.
Since  all points in $[J_nef]$ have the same behavior from $[J_nef]$ to $[J_ne']$, then
their first return time to $[J_n]$ coincide.
\medskip
\ni (ii) By part (i) we only need to prove that 
$V_B^{N_n^{(k-1)}(z)}(z)\in [I_ne_1^{(k)}e_2^{(k)}],
V_B^{N_n^{(k-1)}(y)}(y)\in [I_ne_1^{(k)}e_2^{(k)}]$
for some $e_1^{(k)}e_2^{(k)} \in E(n+1,n+2)$ with $\s(e_1^{(k)})=i^*$. 
The proof is analogous to that of part (i). 
Let us describe the dynamics of a point 
$y\in [I_ne_1...e_{k+1}]$: (1) it moves from  $[I_ne_1...e_{k+1}]$ to the maximal 
path from $v_0$ to $\s(e_{k})$; (2) 
it moves  from this maximal path to the minimal path from $v_0$ to some
$i' \in V_{n+k+1}$, where $i'=\s(succ(e_{k+1}))$ (if $e_k$ is maximal,
then we put $\s(succ(f))=1$ because all minimal edges are connected with $1$);
(3) finally, it moves from this minimal path to the cylinder set 
$[I_ne_1'...e_{k}']$, where $e_1'...e_{k}'$ is a path starting 
at $i^* \in V_n$ and finishing at $i'$. Therefore there is 
$\bar e_1...\bar e_{k} \in E(n+1,n+k)$, $\s(\bar e_1)=i^*$, such that 
$V_B^{N_n^{(1)}(z)}(z) \in [I_n\bar e_1...\bar e_{k}], 
V_B^{N_n^{(1)}(y)}(y) \in [I_n\bar e_1...\bar e_{k}]$.
We conclude by induction.\eop
\centerline{}
\epsfxsize=5truecm
\medskip
\centerline{\epsfbox{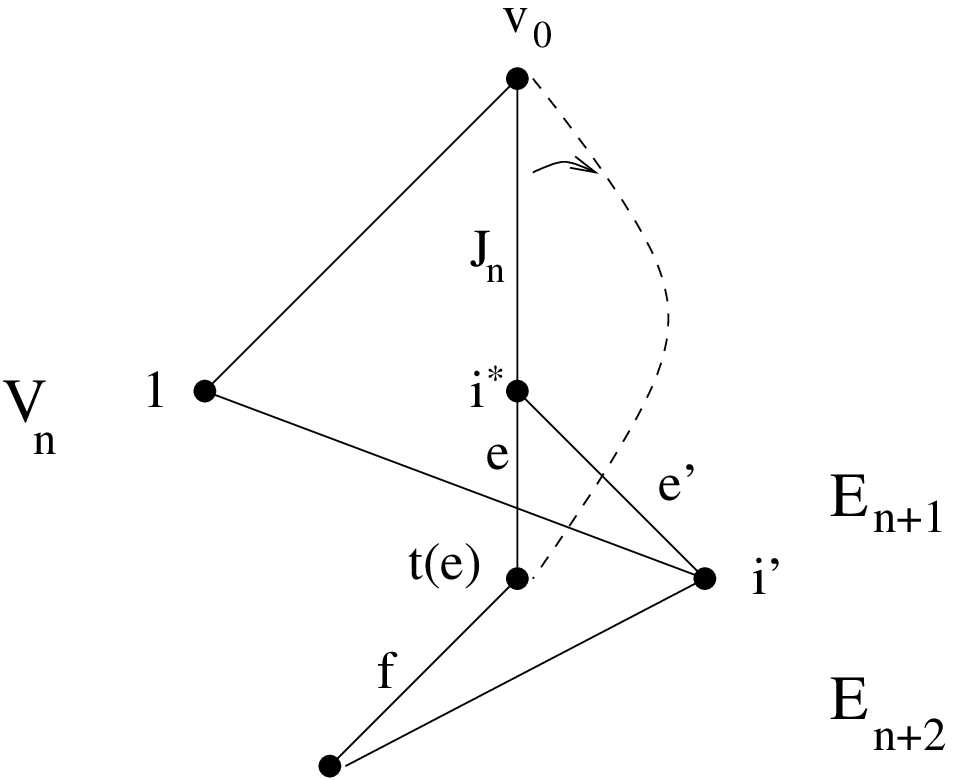}}
\medskip
\centerline{\bf Figure 3}
\medskip
For the sequel, let us fix $n \in {\cal N}$.
Set $\tau_n^{(1)}=\{N_n^{(1)}(x):x \in I_n\}=\{r_1^{(n)},...,r_{l_n}^{(n)}\}$
to be the set of return times to $I_n$
where we are assuming the elements are in increasing order. Also,  
denote $\tau_n^{(1)}(i)=\{x \in I_n: N_n^{(1)}(x)=r_i^{(n)}\}$ for 
$i\in \{1,...,l_n\}$. 
\medskip
\ni{\bf Lemma 2.2} {\it $\mu\{x\in X_B: N_n^{(1)}(x)=k\}=\displaystyle
\sum_{i=1}^{l_n}
1_{\{k\le r_i^{(n)}\}}\cdot \mu( \tau_n^{(1)}(i)).$}
\smallskip
\ni {\bf Proof:} It is clear that 
$${\cal P}_n=\{V_B^{k_i}(\tau_n^{(1)}(i)): 
i\in \{1,...,l_n\}, k_i\in \{0,...,r_i^{(n)}-1\}\}$$ 
is a clopen partition of $X_B$ such that 
$\mu(V_B^{k_i}(\tau_n^{(1)}(i)))=\mu(\tau_n^{(1)}(i))$. It follows that, 
$N_n^{(1)}(x)=k$ if and only if $x \in V_B^{k_i}(\tau_n^{(1)}(i))$
for some $i \in \{1,...,l_n\}$ with $r_i^{(n)}-k_i=k.$ \eop

\ni We denote by $\lfloor \cdot \rfloor$ the integer part 
of a real number.
\medskip
\ni{\bf Lemma 2.3} {\it  For all $t\ge 0$
$$F_n^{(1)}(t)=
\sum_{ef \in E(n+1,n+2):\s(e)=i^*} 
\min \left( {\lfloor {\lambda^{n-1}t \over r(i^*)}\rfloor}, N_n^{(1)}([I_nef]) \right) 
\cdot {r(\t(f)) \over \lambda^{n+1}}.$$}
\smallskip
\ni{\bf Proof:} By Lemma 2.1  the return times to $I_n$ depend
only on the dynamics of points  in the cylinder sets
constructed as the ``continuation''
of $I_n$ by paths of length two.
Hence,  from Lemma 2.2 we get
\medskip
$$\eqalign{
F_n^{(1)}(t) 
&=\sum_{k=1}^{{\lfloor {t \over \mu([I_n]) }\rfloor}} \mu\{x\in X_B: N_n^{(1)}(x)=k\} \cr
&= \sum_{k=1}^{{\lfloor {t \over \mu([I_n]) }\rfloor}} 
\sum_{ef \in E(n+1,n+2):\s(e)=i^*} 1_{\{k\le N_n^{(1)}([I_nef])\}} \cdot \mu([I_nef]) \cr
&= 
\sum_{ef \in E(n+1,n+2):\s(e)=i^*} 
\sum_{k=1}^{{\lfloor {t \over \mu([I_n])}\rfloor}} 
1_{\{k\le N_n^{(1)}([I_nef])\}} \cdot {r(\t(f)) \over \lambda^{n+1}} \cr
&= \sum_{ef \in E(n+1,n+2):\s(e)=i^*}
\min \left( {\lfloor {t \over \mu([I_n]) }\rfloor}, N_n^{(1)}([I_nef]) \right)
\cdot {r(\t(f)) \over \lambda^{n+1}} .}$$
Since $\mu([I_n])={r(i^*) \over \lambda^{n-1}}$,  we conclude the lemma.
\eop
\medskip
Let us compute $N_n^{(1)}([I_nef])$. 
It depends only on the number of times 
the trajectory of a point
in $[I_nef]$  passes through  the minimal path from $v_0$ to a vertex $i \in V_n$ before 
coming  back to $i^*$. 
We call this quantity $c(ef)(i)$.
We remark that this quantity 
does not depend on $n$ because the
diagram is stationary. So $ef \in E(n+1,n+2)$ can be identified 
with some $e'f' \in E(2,3)$.
In addition, when such a trajectory 
passes through  this minimal path then before coming  back to
$i^*$ it has to pass through all paths from $v_0$ to $i$. There are exactly 
$\sum_{k=1}^m M^{n-1}_{k,i}$ of such paths. 
We get, 
$$N_n^{(1)}([I_nef])=\sum_{i=1}^m c(ef)(i) \sum_{k=1}^m M^{n-1}_{k,i}.$$
Let $c(ef)=(c(ef)(i):i\in \{1,...,m\})^T$. In this vector it is ``hidden'' 
the order of the given Bratteli-Vershik system. 
\medskip
We need to compute $\displaystyle \lim_{n\to \infty, n\in {\cal N}}
{N_n^{(1)}([I_nef]) \over \lambda^{n-1}}$ for $ef \in E(n+1,n+2)$. We know from
Perron--Frobenius Theorem (see [HJ]), that 
$\lim_{n\to \infty} {M^{n-1}_{i,j} \over \lambda ^{n-1}}=r(i)l(j)$. Therefore, 
$$\eqalign{ \lim_{n\to \infty, n\in {\cal N}}{N_n^{(1)}([I_nef]) \over \lambda^{n-1}}
&= \lim_{n\to \infty, n\in {\cal N}} \sum_{i=1}^m c(ef)(i) \sum_{k=1}^m M_{k,i}^{n-1} \cdot 
{1 \over \lambda^{n-1}}\cr
                                &= \sum_{i=1}^m c(ef)(i) \sum_{k=1}^m r(k)l(i)  = \sum_{i=1}^m c(ef)(i) l(i)= \bar c(ef). \cr}
$$
\medskip
Let  $L=|\{\c(ef):ef \in E(n+1,n+2)\}|$. Since the Bratteli diagram
is stationary $L=|\{\c(ef):ef \in E(2,3)\}|$.
We set $\{\c(ef):ef \in E(2,3)\}=\{\c(e_if_i): i\in \{1,...,L\}\}$.
We also assume that
$0 < c_1=\c(e_1f_1)<...<c_{L}=\c(e_{L}f_{L})$.
We set
$S(i)=\{ef \in E(n+1,n+2) : \s(e)=i^*, \c(ef)=c_i \}$ for each $i\in  \{1,...,L\}$.
Paths $e_if_i$ and sets $S(i)$
can be assumed to be the same for every $n\ge 1$ since
the Bratteli diagram is stationary. Finally put $d_0=0,d_1=c_1\cdot r(i^*),...,d_{L}= c_{L} 
\cdot r(i^*)$, and $d_{L+1}=\infty$.
\bigskip				
\ni {\bf Theorem 2.4} {\it Let  
$d_j\le t < d_{j+1}$ and $j\in \{0,...,L\}$. Then,  
the limit laws are the piecewise 
linear functions given by Figure 4 which can be described as follows, 
$$F^{(1)}(t)=\lim_{n\to \infty, n\in {\cal N}} F_n^{(1)}(t)= 
\sum_{ef \in \cup_{i=1}^{j}S(i)} \c(ef)  {r(\t(f)) \over \lambda^2 } + 
{t \over r(i^*) }
\cdot  \sum_{ef \in \cup_{i=j+1}^L S(i)} {r(\t(f)) \over \lambda^2 }$$
and
$$F^{(k)}(t)=\lim_{n\to \infty, n\in {\cal N}} F_n^{(k)}(t)=
\sum_{e_1...e_{k+1} \in E(2,k+2): \s(e_1)=i^*,e_1^{(k)}e_2^{(k)}\in \cup_{i=1}^jS(i)} 
\c(e_1e_2)  {r(\t(e_{k+1})) \over \lambda^{k+1}}.$$
The convergence is also uniform in any closed interval 
$I \subseteq ]d_j,d_{j+1}[$.}
\medskip
\epsfxsize=7truecm
\medskip
\hskip 0.7cm
\epsfbox{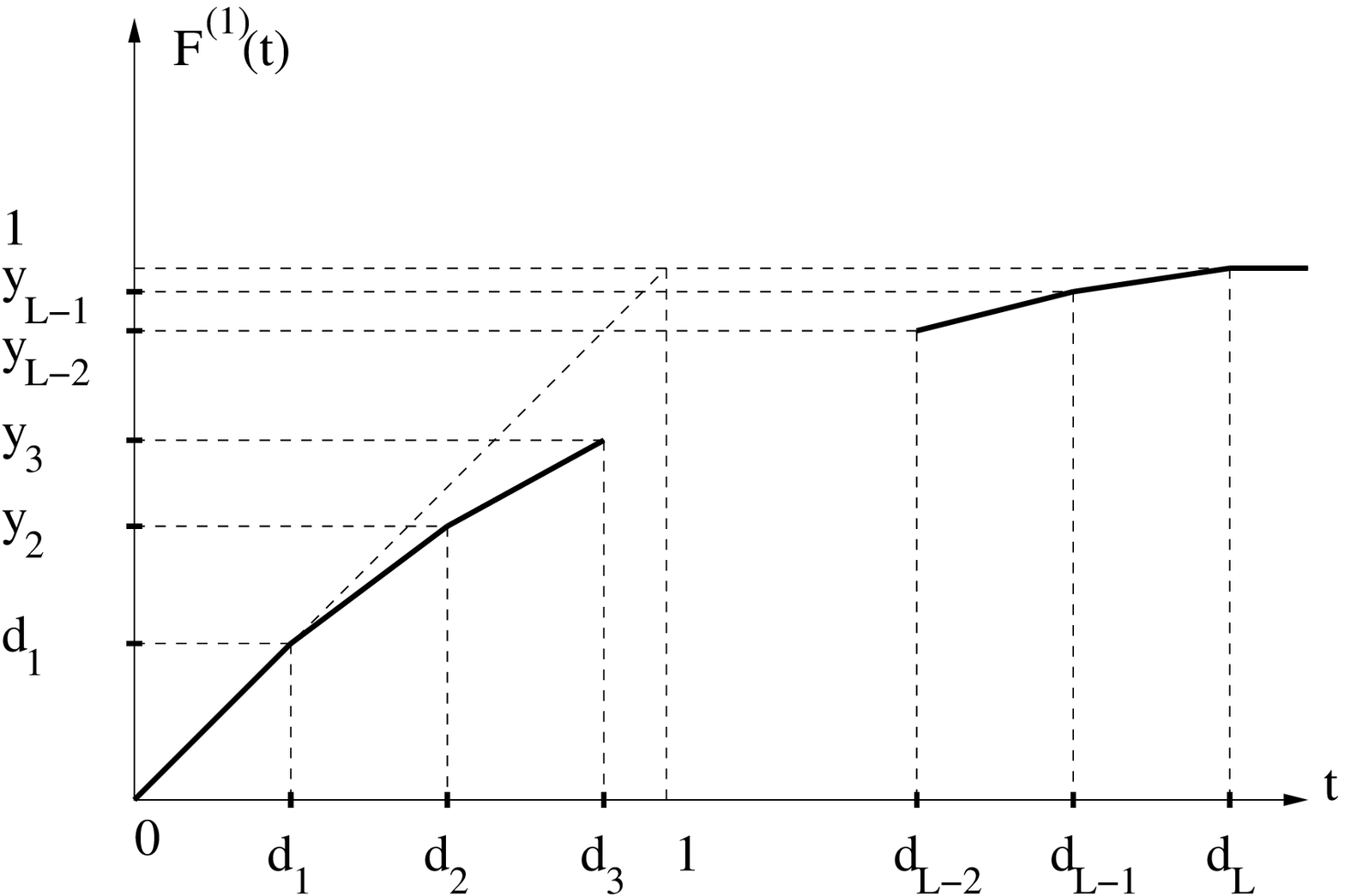}
\hskip 1cm
\epsfxsize=7truecm
\epsfbox{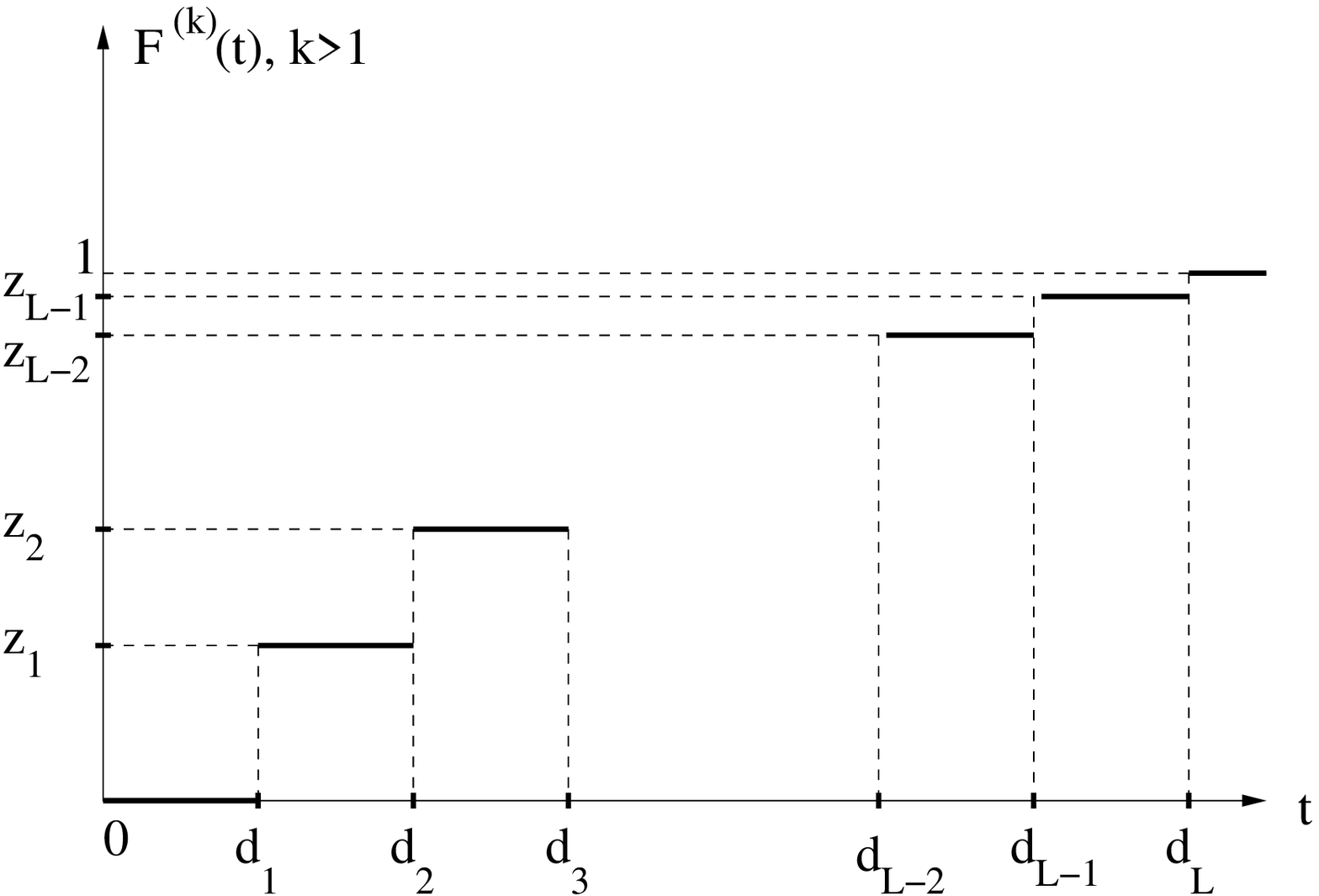}

\medskip
\centerline{\bf Figure 4}
\medskip
\ni {\bf Proof:} 
We start with the computation of the limit law for the first entrance time. 
Fix $d_j\le t < d_{j+1}$ where $j\in \{0,...,L\}$.
From Lemma 2.3 we get 
$$ F_n^{(1)}(t) = \sum_{ef \in E(n+1,n+2):\s(e)=i^*}
\min \left( {1 \over \lambda^{n-1}} {\lfloor {\lambda^{n-1}t \over r(i^*)}\rfloor}, 
{ N_n^{(1)}([I_nef]) \over \lambda^{n-1}}  \right)
\cdot {r(\t(f)) \over \lambda^{2}}$$
Since $E(n+1,n+2)$ can be identified with $E(2,3)$, taking limit in $n \in {\cal N}$ we conclude, 
$$\eqalign{
F^{(1)}(t)&=\lim_{n\to \infty, n \in {\cal N}} F_n^{(1)}(t)= 
\sum_{ef \in E(2,3):\s(e)=i^*} \min \left( {t \over r(i^*)},\c(ef) \right) \cdot 
{r(\t(f)) \over \lambda^2} \cr
&=\sum_{ef \in \cup_{i=1}^j S(i)} \c(ef) \cdot {r(\t(f)) \over \lambda^{2}} +
{t\over r(i^*)}
\sum_{ef \in \cup_{i=j+1}^{\ell} S(i)} {r(\t(f)) \over \lambda^{2}}.}$$
\medskip
Now we compute the limit for $F_n^{(k)}(t)$, $k\ge 2$. That is, 
$$F_n^{(k)}(t)=\sum_{s=1}^{\lfloor {\lambda^{n-1} t \over r(i^*)} \rfloor}
\mu\{x\in X_B: N_n^{(k)}(x)-N_n^{(k-1)}(x)=s\}.$$
By Lemma 2.1, the difference $N_n^{(k)}(x)-N_n^{(k-1)}(x)$ only depends on the 
cylinder set $[I_ne_1...e_{k+1}]$ which contains $x$. Indeed, if 
$x,y \in [I_ne_1...e_{k+1}]$ then $N_n^{(1)}(V_B^{N_n^{(k-1)}(x)}(x))=N_n^{(1)}(V_B^{N_n^{(k-1)}(y)}(y))$,
because $V_B^{N_n^{(k-1)}(x)}(x),V_B^{N_n^{(k-1)}(y)}(y)\in [I_ne_1^{(k)}e_2^{(k)}]$.
Then, 
$$F_n^{(k)}(t)=\sum_{e_1...e_{k+1}\in E(n+1,n+k+1): \s(e_1)=i^*} \hskip -0.5truecm
1_{\{N_n^{(1)}([I_ne_1^{(k)}e_2^{(k)}])
\le {\lfloor {\lambda^{n-1} t \over r(i^*)} \rfloor}\}}  
{N_n^{(1)}([I_ne_1e_2])  
r(\t(e_{k+1})) \over \lambda^{n+k}}.$$
Take $d_j\le t < d_{j+1}$, $j\in\{0,...,L\}$. In a similar way as we did
for $F_n^{(1)}(t)$ we get that, 
$$F^{(k)}(t)=\lim_{n\to \infty, n\in {\cal N}} F_n^{(k)}(t)=\sum_{e_1...e_{k+1}\in E(2,k+2):\s(e_1)=i^*,
e_1^{(k)}e_2^{(k)} \in \cup_{i=1}^j S(i)} \c(e_1e_2) \  {r(\t(e_{k+1})) \over \lambda^{k+1}}.$$

Since the number of return times is bounded, 
a standard compactness argument proves that the convergences are also uniform. \eop 
\medskip
Let us point out that the limit laws  provided in the theorem do not depend on 
the explicit sequence of cylinder sets $(I_n:n \in \NN)$ considered, but only on the vertex $i^*$.
Then each $i \in \{1,..,m\}$ defines its own family of  limit laws. Consequently
for $k\ge 1$, $(F_n^{(k)}: n\in {\cal N})$ converges if and only if 
the limit laws defined by 
the terminal vertices of $I_n$, for all $n \in {\cal N}$  
large enough, coincide.
\medskip
From the computation in the proof of Theorem 2.4 and some considerations from matrix theory 
(see [HJ] Theorem 8.5.1) we get the 
following convergence rate.
\medskip
\ni{\bf Corollary 2.5} 
{\it There exist positive constants $\gamma, C, D$  such that  $\gamma < \lambda$ and
$$\sup_{t \in \RR} \left | F_n^{(1)}(t) - F^{(1)}(t) \right | \leq C \left ( {\gamma \over \lambda} \right )^n 
\hbox{ and } \sup_{t \in \RR} \left | F_n^{(k)}(t) - F^{(k)}(t) \right | \leq D  \left 
( {\gamma \over \lambda} \right )^n. $$ }

\medskip
\ni{\bf Example 2 (left to right order):} In this example we consider 
Cantor minimal systems given by stationary 
Bratteli-Vershik diagrams satisfying conditions 
(H1),(H2),(H3) and with increasing  order from left to right. That is, for any 
$n\ge 1$ and for any 
$e,f \in E_n$, if $\t(e) \leq  \t(f)$ then  $\s(e)\leq \s(f)$, 
where in all the $V_n$ we put the natural order of $\{1,..,m\}$ (see Figure 5).
\medskip
\epsfxsize=6truecm
\medskip
\centerline{\epsfbox{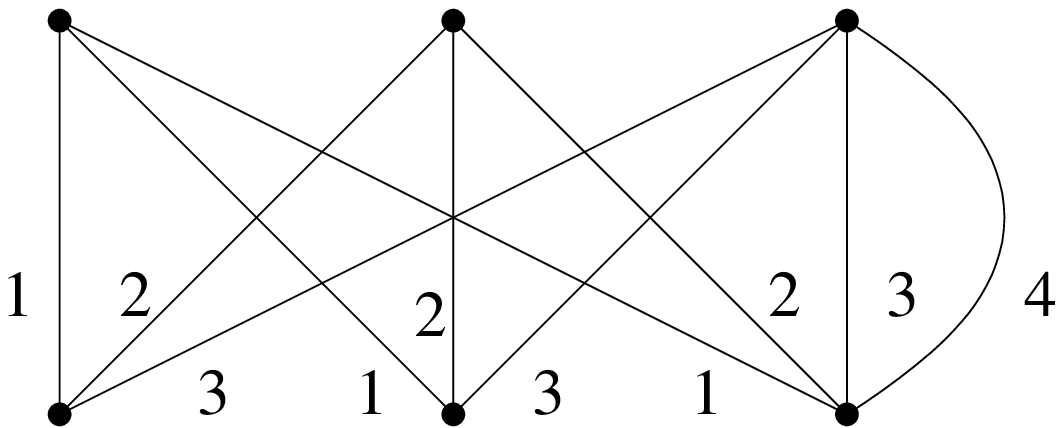}}
\medskip
\centerline{\bf Figure 5}
\medskip

Put $i^*=1$ and let $I_n$ be the minimal path 
from $v_0$ to $1\in V_n$. Then ${\cal N}=\NN$. It is not difficult to see
that return times to $I_n$ are constants over $[I_ne]$ for 
$e\in E_{n+1}$, $\s(e)=1$. Let us fix one of such $e$ and put $j=\t(e)$. 
If $e$ is not a maximal edge with respect to 
the set $G(j,n+1)$ of edges in $E_{n+1}$ from $1 \in V_n$ to $j \in V_{n+1}$, 
then $N_n^{(1)}([I_ne])=\sum_{k=1}^m M_{k,1}^{n-1}$, and if
$e$ is a maximal edge, with respect to the same set of edges, 
we get $N_n^{(1)}([I_ne])=\sum_{k=1}^m M_{k,j}^{n}+(1-M_{1,j})
\sum_{k=1}^m M_{k,1}^{n-1}$.
Dividing by $\lambda^{n-1}$ and taking the limit when $n$ tends to infinity 
we obtain 
$$\{d_1,...,d_{L}\}=\{ r(1)l(1), r(1) l(1) + r(1) \left ( \lambda l(j) - M_{1,j} l(1) ) \right ), j\in \{1,...,m\}\}.$$
Also,  $\c(e)=l(1)$ if $e$ is not maximal in $G(j,n+1)$  and 
$\c(e)=  l(1)+\lambda l(j)- M_{1,j}l(1)$ if it is maximal. Then the limit laws can be deduced from
the general statement in Theorem 2.4.

\bigskip
\bigskip
\ni{\bf 3. Limit laws for non stationary Bratteli--Vershik systems.}
\bigskip
\bigskip
In this section we will compute the limit laws for some 
minimal Cantor 
systems given by non stationary Bratteli diagrams. First 
we give a general formula and then we apply it to  linearly reccurent 
subshifts, odometers and sturmian subshifts. 
As in the  stationary case 
we will fix some properties of the ordered Bratteli diagrams. 
For any Cantor minimal system these properties hold after contracting and microscoping levels
of a given Bratteli-Vershik representation of the system. 
\medskip
Let $(X_B,V_B)$ be a Bratteli-Vershik system and fix a $V_B$-invariant probability measure
$\mu$, where 
$\displaystyle {B}=(\cup_{i\ge 0} V_i, \cup_{i\ge 1} E_i, \preceq)$ with
$V_i=\{1,...,m_i\}$, $i\ge 1$, $V_0=\{v_0\}$. Recall that $(M^{(k)}: k\ge 1)$
are the incidence matrices of levels. Furthermore the 
following properties hold:
\medskip
(H1) the incidence matrices $(M^{(k)}: k\ge 1)$  of $B$ has strictly positive coefficients;
\smallskip
(H2) for every vertex $i\in V_1$ there is a unique edge from $v_0$ to $i$;
\smallskip
(H3) $\forall i \ge 1$, $\forall j \in \{1,...,m_{i}\}$, 
$e=\min \{f \in E_i:\t(f)=j \} \Rightarrow \s(e)=1 \in V_{i-1}$.
\medskip
These conditions allow to prove a version of Lemma 2.1 for general  minimal Cantor systems.
The proof is left to the reader. 
\medskip
\ni{\bf Lemma 3.1} {\it Let $I=[x_1,...,x_n]$ be a cylinder set in $(X_B,V_B)$ with $i^*=\t(x_n)$.

\ni (i) Let $ef\in E(n+1,n+2)$ with $\s(e)=i^* \in V_n$.
Then, for any $z,y \in [Ief]$ we have  $N_I^{(1)}(z)=N_I^{(1)}(y)$.
\smallskip
\ni (ii) Let $e_1...e_{k+1} \in E(n+1,n+k+1)$ with $\s(e_1)=i^*$. Then, there is
$e_1^{(k)}e_2^{(k)} \in E(n+1,n+2)$ with $\s(e_1^{(k)})=i^*$  such that for any
$z,y \in [Ie_1...e_{k+1}]$, $N_I^{(k)}(z)=N_I^{(k)}(y)$ and
$V_B^{N_I^{(k-1)}(z)}(z),V_B^{N_I^{(k-1)}(y)}(y)\in [Ie_1^{(k)}e_2^{(k)}]$.}
\medskip
The last lemma and similar considerations as those made in 
the previous section imply that,
$$F_I^{(1)}(t)=\sum_{ef\in E(n+1,n+2):\s(e)=i^*} 
\min \left ( {\lfloor {t \over \mu(I) }\rfloor}, N_I^{(1)}([Ief]) \right) 
\cdot \mu([Ief])$$
and
$$
F_I^{(k)}(t)=\hskip -6pt \sum_{e_1...e_{k+1} \in E(n+1,n+k+1):\s(e_1)=i^*} \hskip -6pt
1_{\{N_I^{(1)}([Ie_1^{(k)}e_2^{(k)}]) \le {\lfloor {t \over \mu(I) }\rfloor}\}}
N_I^{(1)}([Ie_1e_2]) \cdot \mu([Ie_1...e_{k+1}]). 
$$
\medskip

\ni{\bf Example 3 (Linearly recurrent subshifts):}   
An example  of non stationary Bratteli-Vershik systems 
are linearly recurrent subshifts introduced in [D]. 
They  can be represented by   ordered Bratteli diagrams verifying 
conditions (H1), (H2), (H3),  such that 
for all $n \in \NN$, $|V_n|=|V_{n+1}|$ and  $|E_n|\le K$ where $K$ is  a universal constant.
In addition, it can be proved (following  the same lines in [D])
that  there is a constant $\bar K$ such that for every cylinder set $I$ and 
every $x \in X_B$, $\mu(I) N_I^{(1)}(x) \le \bar K$.
Therefore, once we  fix $k\ge 1$,  we can get a subsequence 
${\cal N}_k \subseteq \NN$ for which  
$\lim_{n\to \infty, n\in {\cal N}_k} F ^{(k)}_n$
exist and is a piecewise linear function as in the case 
of Theorem 2.4. 
\medskip
The following example is neither stationary nor linearly recurrent.
\bigskip
\ni{\bf Example 4 (Odometer):} Let $(X,T)$  
be the odometer with base $(p_n:n\in \NN)$. 
The level $n$ of the classical Bratteli-Vershik representation of odometers is given in Figure 6(a).
\medskip
\epsfxsize=10truecm
\medskip
\centerline{\epsfbox{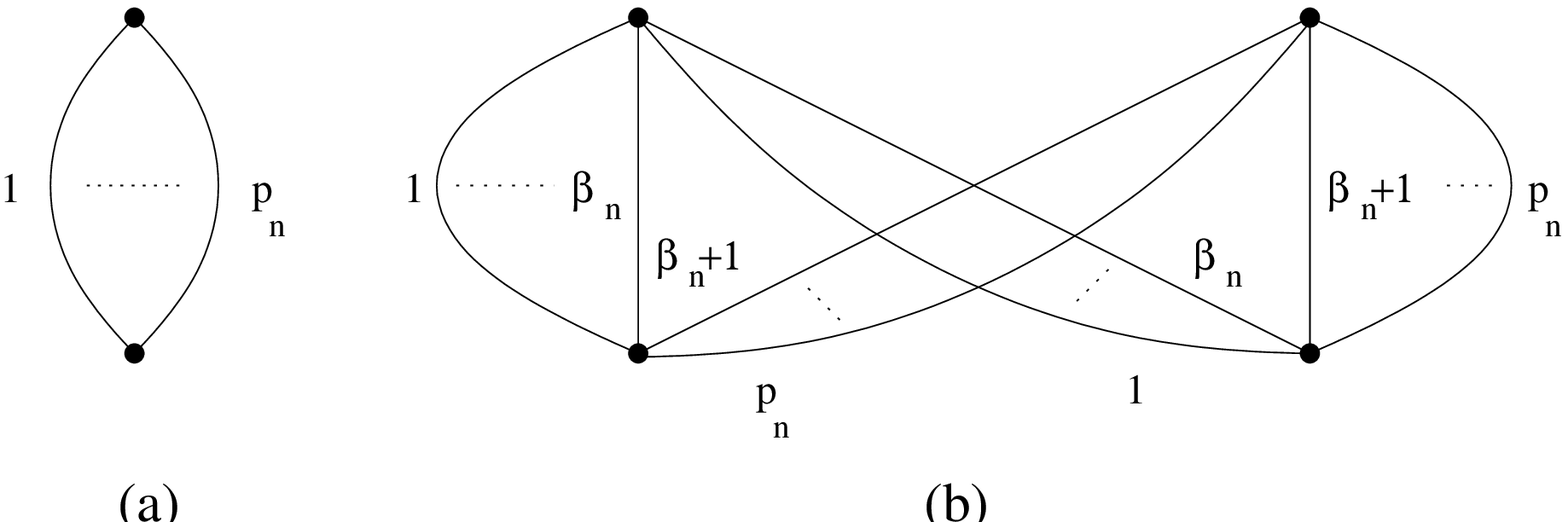}}
\medskip
\centerline{\bf Figure 6}
\medskip
In this case, if we take $I_n$ to be any cylinder set of length $n$, then 
the limit law of first entrance time is a uniform law in $[0,1]$ and for the 
$m$-th return time it is a discrete distribution concentrated in $1$.
\medskip
Let $\beta \in (0,1)$.
Another representation by means of Bratteli diagrams 
of an odometer is given by Figure 6(b). 
We set $\beta_n=\lfloor \beta  p_n \rfloor$.
Let $(I_n:n \in \NN)$ be a sequence of cylinders set induced by 
a point $x^* \in X_B$.
The unique ergodic measure of the system is given by
$\mu(I_n)={\beta_n  \over q_{n}}$, where $q_n=p_1\cdot...\cdot p_n$. 
There are two values for the return times to $I_n$: 
$q_{n-1}$ and $q_{n-1} +q_{n-1} (p_n - \beta_n)$. Then, 
$d_1^{(n)}= \mu(I_n) q_{n-1}= {\beta_n  \over p_n}$ and 
$d_2^{(n)}= \mu(I_n) \left ( q_{n-1} +q_{n-1} (p_n - \beta_n) \right )= {\beta_n \over p_n} (1+p_n-\beta_n)$.
If there is a subsequence $(n_i:i\in \NN)$ such that $p_{n_i}=p$ then 
$d_1=d_1^{(n_i)}={\lfloor \beta p \rfloor \over p}$ and 
$d_2=d_2^{(n_i)}={\lfloor \beta p \rfloor \over p} (1+p-\lfloor \beta p \rfloor)$.
In this case the limit law for the first entrance time is given by the piecewise 
linear function in Figure 7(a). This limit is uniform.
If $\lim_{n\to \infty}p_n=\infty$ then  
$\lim_{n\to \infty}d_1^{(n)}=\beta$ and $\lim_{n\to \infty}d_2^{(n)}=\infty$. 
Consequently the pointwise limit is given by Figure 7(b) and the limit is not uniform.
\medskip
\epsfxsize=13truecm
\medskip
\centerline{\epsfbox{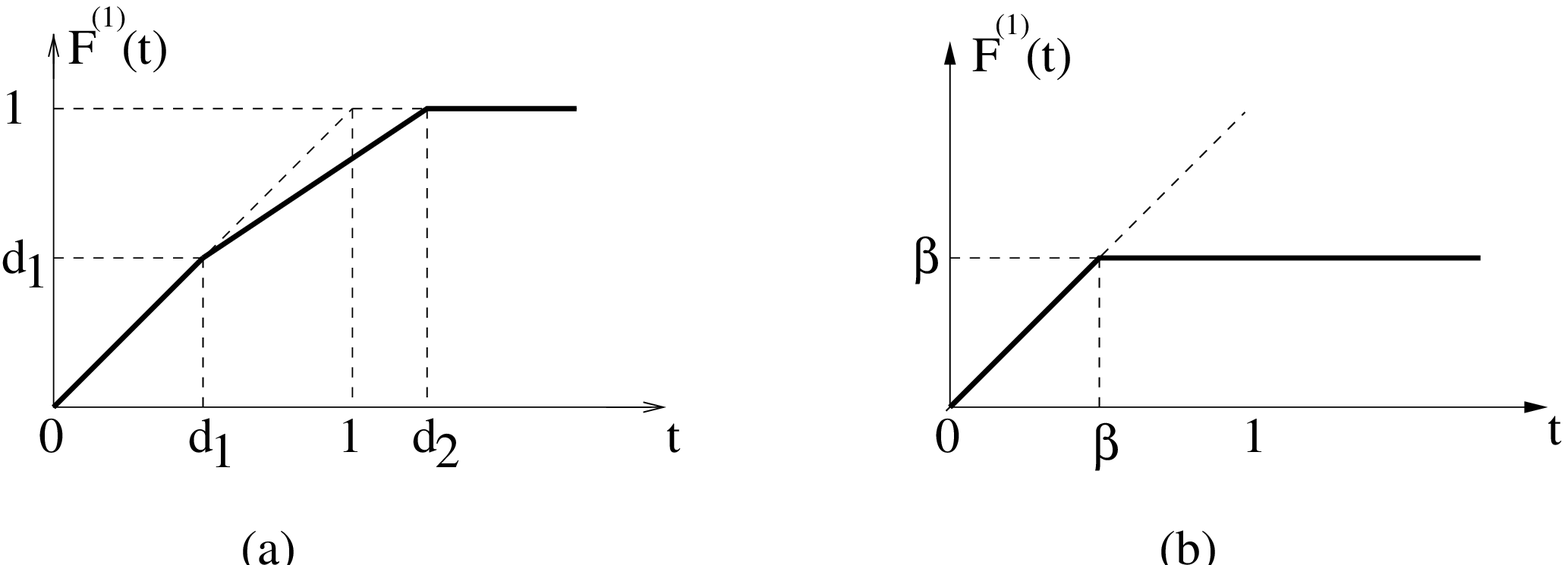}}
\medskip
\centerline{\bf Figure 7}
\bigskip
\ni{\bf Example 5 (Sturmian subshifts):} 
This example is motivated by the results in [CdF], where the authors computed 
the limit laws of entrance times for rotations of the circle. In the context of 
subshifts they correspond to Sturmian systems. Not surprisingly the results we obtain 
here are analogous.
\medskip
Let $0<\alpha <1$ be an irrational number. We define the map 
$R_{\alpha }:\left[ 0,1 \right[ \rightarrow \left[ 0,1 \right[ $ 
by 
$R_{\alpha} \left ( t \right ) = t + \alpha$ 
(mod 1) and the map $I_{\alpha }:\left[ 0,1\right[ \rightarrow \left\{ 0,1\right\}$
by $I_{\alpha }\left( t\right) =0$ if $t\in \left[ 0,1-\alpha \right[$
and $I_{\alpha }\left( t\right) =1$ 
otherwise. 
Let 
$\Omega _{\alpha }=\overline{\left\{ \left. \left( I_{\alpha }\left( R^{n}_{\alpha }\left( t\right) 
\right) \right) _{n\in \ZZ }\right| t\in \left[ 0,1\right[ \right\} }
\subset \left\{ 0,1\right\} ^{\ZZ }$.
The subshift $\left( \Omega _{\alpha },\sigma \right)$ is called  Sturmian
subshift (generated by $\alpha$) and its elements are called Sturmian sequences.
There exists a factor map (see [HM]) 
$\gamma :\left( \Omega _{\alpha },\sigma \right) \rightarrow \left( \left[ 0,1\right[ ,R_{\alpha }\right)$  
such that,

(1)  $\left| \gamma ^{-1}\left( \left\{ \beta \right\} \right) \right| =2$  if
$\beta \in \left\{ \left. n\alpha \right| n\in \ZZ \right\}$
and

(2) $\left| \gamma ^{-1}\left( \left\{ \beta \right\} \right) \right| =1$ otherwise.

This map induces a measure-theoretical isomorphism.
It is also well known that Sturmian systems are uniquely ergodic 
and the symbolic complexity is $n+1$ [HM].
\medskip
A  Bratteli-Vershik representation
for Sturmian subshifts is presented in [DDM]. 
It works as follows. There is a sequence of positive 
integers $(d_k:k\ge 1)$ such that level $k$ of the Bratteli diagram 
is given by either block (a) or  block (b) of Figure 8. In addition, it does not 
exist two consecutive levels ordered like block (a) in Figure 8.
We notice that blocks (a) and (b) have the same 
incidence matrix $M^{(k)}=N_{d_k}=
\left[\matrix{ d_k & 1 \cr 1 & 0 \cr } \right ]$ but different orders.
Also, the continued fraction expansion of $\alpha$ and $\beta=[0:d_1,d_2,\dots]$
are eventually equal.  For the sequel we fix such representation  where  we identify
the set of vertices $V_n$ with  $\{1,2\}$.
\medskip
\epsfxsize=4truecm
\medskip
\hskip 2cm
\epsfbox{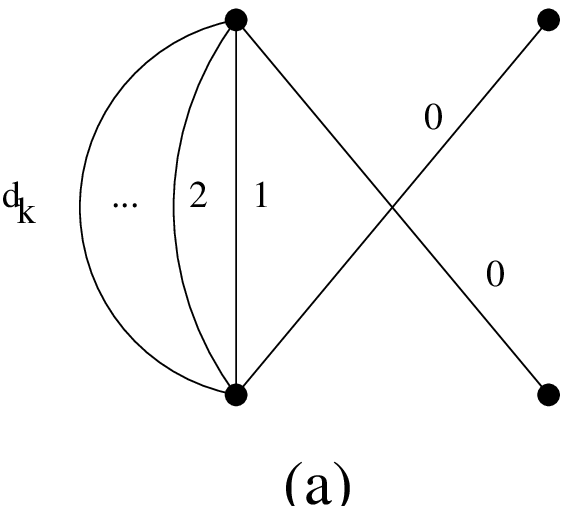}
\hskip 1cm
\epsfxsize=4truecm
\epsfbox{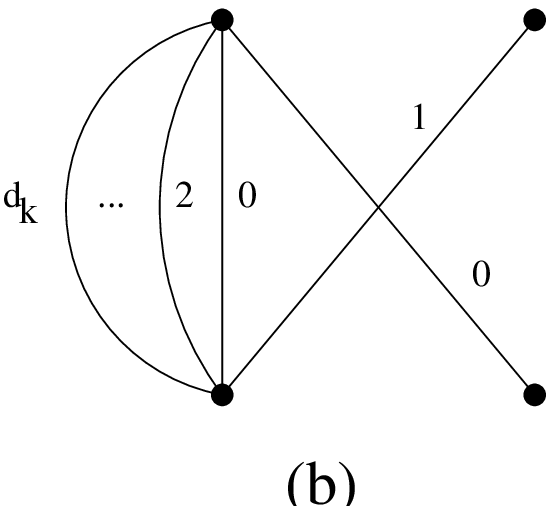}
\medskip
\centerline{\bf Figure 8}
\medskip
We compute the limit laws of the first entrance time 
for a sequence of cylinders $I_n=[x_1...x_n]$ such that 
infinitely many times $\t(x_n)=1$. 
Let ${\cal N}=\{n\in \NN : \t(x_n)=1 \}$.
Denote 
$F_{I_n}^{(1)}(t)=F_n^{(1)}$ and $N_{I_n}^{(1)}(t)=N_n^{(k)}$.
Let us notice that the ordered Bratteli diagram constructed in [DDM] does not 
verify properties (H1), (H2) and (H3) but the results of Lemma 3.1 
still hold. Then to compute limit laws we need to know:
$N_n^{(1)}([I_nef])$ for every $ef \in E(n+1,n+2)$ such that 
$\s(e)=1$ and $\mu([I_n])$, where $\mu$ is the unique invariant measure of the 
Sturmian subshift. In that purpose we need to recall some results of continued 
fraction theory. First, 
$N_{d_1} \cdot ... \cdot N_{d_{k-1}}=\left [ 
\matrix{ p_{k-1} & p_{k-2} \cr q_{k-1} & q_{k-2} } \right ]$ where 
${p_k \over q_k}=[0:d_1,...,d_k]$ is the classical  
approximation of $\beta$ (see [HW]). From this expression we deduce that 
$\mu([y_1...y_k])={1 \over \beta +1} |\beta q_{k-2} - p_{k-2}|$ if 
$\t(y_k)=1$ and $\mu([y_1...y_k])={1 \over \beta +1} |\beta q_{k-1} - p_{k-1}|$ if
$\t(y_k)=2$,  and that 
$q_k|\beta q_k - p_k|={G^k(\beta) \over 1+ (q_{k-1} / q_k) G^k(\beta)}$ where 
$G$ is the Gauss map: $G(\beta)=\left \{ 1/ \beta \right \}$ 
(the fractional part of $1 / \beta$).
\medskip
Let us fix $n \in {\cal N}$. Looking at the diagram we verify that 
there are two possible values for $N_n^{(1)}([I_nef])$ with $ef \in E(n+1,n+2)$ and 
$\s(e)=1$: $p_{k-1}+q_{k-1}$ and $p_{k-1}+q_{k-1}+p_{k-2}+q_{k-2}$. Suppose 
there is a subsequence ${\cal M} \subseteq {\cal N}$ such that 
$$ \lim_{n\to \infty, n\in {\cal M}} {q_{n-3} \over q_{n-2}}=w \hbox{ and } 
   \lim_{n\to \infty, n\in {\cal M}} G^{n-2}(\beta)=\theta.$$
These 
conditions are exactly the same as those provided in [CdF] to have a non trivial
limit. We get
$$
h_1=\lim_{n\to \infty, n\in {\cal M}} (p_{n-1}+q_{n-1})\cdot \mu([I_n])=
{ \lfloor{ 1\over \theta} \rfloor \theta \over 1+ \theta w} \quad  \hbox{ and }$$ 
$$h_2=\lim_{n\to \infty, n\in {\cal M}} (p_{n-1}+q_{n-1}+p_{n-2}+q_{n-2})\cdot \mu([I_n])=
{ (1 + \lfloor{ 1\over \theta} \rfloor) \theta \over 1+ \theta w}.
$$
Then, the limit $F^{(1)}(t)=\lim_{n\to \infty, n\in {\cal M}} F^{(1)}_n(t)$ is the continuous 
piecewise linear function given by Figure 9. An analogous computation yields  
$F^{(k)}(t)$.
\medskip
\epsfxsize=5truecm
\medskip
\centerline{\epsfbox{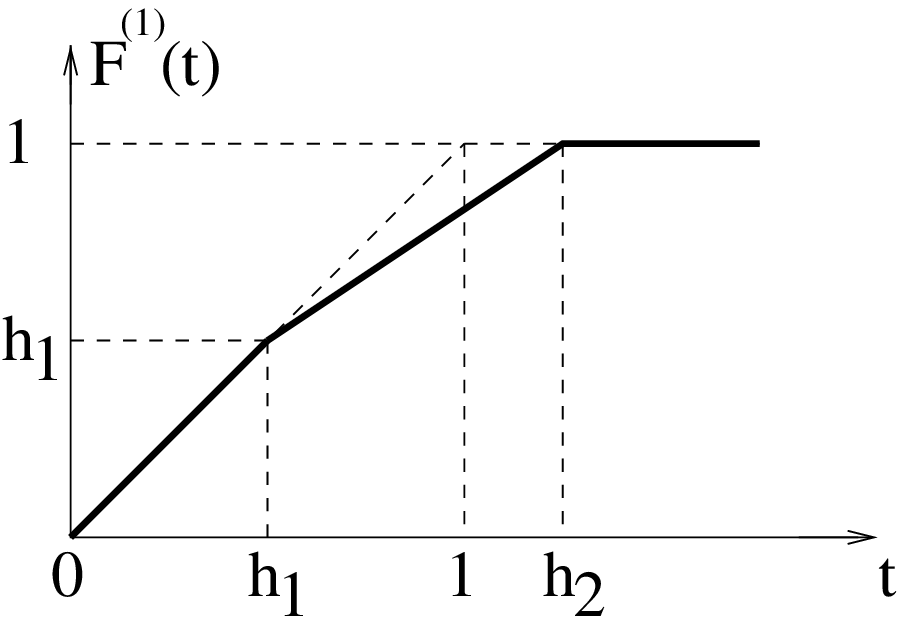}}
\medskip
\centerline{\bf Figure 9}
\vfill\eject
\bigskip
\bigskip
\ni{\bf 4. Point process induced by entrance times.}
\bigskip
\bigskip
Let $(X,T)$ be a minimal Cantor system 
and $\mu$ a $T$-invariant probability 
measure. Consider a decreasing family of clopen sets
of $X$, $(I_n:n \in \NN)$ such that $\cap_{n\in \NN} I_n=\{x^*\}$.
For each $x \in X$ and $k\ge 2$ define 
$T_n^{(k)}(x)=N_n^{(k)}(x)-N_n^{(k-1)}(x)$, $T_n^{(1)}(x)=N_n^{(1)}(x)$.
Denote by $\delta_t$ the Dirac measure at the point $t \in \RR$.
The point process  $\tau_n:X\to {\cal M}[0,\infty)$ defined by this sequence of renewal times is 
$$\tau_n(x)=\sum_{k\ge 1} \delta_{N_n^{(k)}(x) \mu(I_n)},$$
where $x$ is randomly chosen with respect to $\mu$
and ${\cal M}[0,+\infty)$ is the set of 
$\sigma$-finite measures on $[0,+\infty)$.
In this section we consider the problem whether this point process converges in law. 
To see that we have  to compute the limit of 
$$F_{1,..,p}^{(n)}(t_1,...,t_p)=\mu\{x\in X: T_n^{(1)}(x) \mu(I_n)\le t_1,...,T_n^{(p)}(x) \mu(I_n) \le t_p\}$$
when $n$ tends to infinity in some subsequence ${\cal N}\subseteq \NN$, for all 
$p\in \NN$ and for all $(t_1,...,t_p) \in \RR^p$ (see [N]). 
\medskip
We will focus on the stationary case. We set the same notations used 
in section 2. Following  the same lines of the proof of Theorem 2.4 we get, 

$$\eqalign{
F^{(n)}_{1,...,p} (t_1,...,t_p) =&\sum_{e_1...e_{p+1} \in E(n+1,n+p+1):\s(e_1)=i^*}
\min \left( \lfloor { t_1 \lambda^{n-1} \over r(i^*) } \rfloor,  N_n^{(1)}([I_ne_1e_2])\right ) \cr
&\prod_{k=2}^p 1_{\{ N_n^{(1)}([I_ne_1^{(k)}e_2^{(k)}]) < {t_k \lambda^{n-1} \over r(i^*)  } \}}
{r(\t(e_{p+1})) \over \lambda^{n+p+2}} }
$$
where for $k \in \{2..,p\}$, $e_1^{(k)}e_2^{(k)} \in E(n+1,n+2)$ is the unique path such that 
$T^{N_n^{(k-1)}(x)}(x)$ belongs to   
$[I_ne_1^{(k)}e_2^{(k)}]$ for every $x \in [I_ne_1...e_{p+1}]$.
Therefore the point process $\tau_n$ converges to (a priori) a non-stationary point process
para\-metri\-zed by the Bratteli-Vershik  diagram with distribution
$$\eqalign{
F_{1,...,p} (t_1,...,t_p) =&\sum_{e_1...e_{p+1} \in E(1,p+1):\s(e_1)=i^*}
\min \left( t_1 r(i^*), \c(e_1e_2) \right ) \cr
&\prod_{k=2}^p 1_{ \{ \c(e_1^{(k)}e_2^{(k)}) \le {t_k\over r(i^*)  }\}}
{r(\t(e_{p+1})) \over \lambda^{p+3}}. }
$$
\bigskip
\bigskip
\ni{\bf 5. Final comments and questions.}
\bigskip
\bigskip
The results presented in this paper together with those obtained 
in [CdF], and the results in [CC], [CG], [HSV] and [P] (among others), 
show two extreme behaviors for some  families
of sequences $(I_n:n\in \NN)$. 
In the first case the limit laws    
are piecewise linear functions 
and in the others they are exponential laws (Poisson laws). 
\medskip
All subshifts considered in this paper, 
substitutions subshifts, linearly recurrent subshifts 
and Sturmian subshifts, 
share  at least two common features that are intimately 
related to the existence of a 
linear limit law for a given sequence of cylinder sets
$(I_n:n\in \NN)$: their ordered Bratteli diagrams 
are ``universally bounded'' and 
their symbolic complexity is sub-linear ([Q], [DHS], [HM]). The second condition is  
behind the fact that $\mu(I)N_I^{(1)}(x)$ is 
bounded independently of $I$ and $x$.
A natural question is whether piecewise linear limit 
laws are characteristic of systems verifying such conditions. 
\medskip
The result of Lacroix [L] about limit laws 
for the first return time tells us that any distribution 
function can be obtained as a limit law if we choose correctly
the sequence $(I_n:n\in \NN)$. 
On the other hand, it seems that the community working on this topic agrees 
that exponential limit laws should be  characteristic of 
mixing systems with positive entropy 
and piecewise linear limit laws should be characteristic of 
subshifts with sublinear complexity. 
The result of Lacroix shows that 
these facts need to be clarified. One possible direction is to explore 
which are the ``natural'' sequences $(I_n:n\in \NN)$. 
Once these sequences are defined we could  ask whether 
there is a dynamical system which  
limit laws are in between piecewise linear functions and 
exponential laws.
A natural class to consider is that of Toeplitz
subshifts. 
They can have positive or  zero entropy [W] and 
they can be represented by means of ordered Bratteli diagrams 
with a nice structure [GJ].
For example, which limit laws can we obtain for  Toeplitz systems with polynomial symbolic complexity 
(see [CK])  ?
\bigskip
\bigskip
\ni {\bf Acknowledgements:}  this paper 
was partly written during the first author visit to the 
Centro de Modelamiento Matem\'atico of the Universidad de Chile and 
the second author visit to the LAMFA-CNRS of 
the Universit\'e de Picardie Jules Verne. 
The support and kind hospitality
of both institutions were very much appreciated. 
The authors had partial support from ECOS-Conicyt program, 
C\'atedra Presidencial 
fellowship and Fondap program in Applied Mathematics.
\vfill\eject
\bigskip
\bigskip
\ni{\bf References.}
\bigskip
\bigskip

\item{[BJKR]} O. Bratteli, P.E.T. Jorgensen, K.H. Kim, F. Roush, 
Non-stationarity of isomorphism between AF algebras defined by stationary Bratteli diagrams, 
preprint (1999).

\item{[CK]}  J. Cassaigne, J. Karhum\"aki, Toeplitz words, generalized periodicity and periodically iterated
morphisms, European J. Combin. 18 (1997), no. 5, 497--510.

\item{[C]} Z. Coelho, Asymptotic laws for symbolic dynamical systems, 
London Mathematical Society, Lecture Notes Series 279, Cambridge University Press (2000), 
123--165.

\item{[CC1]} Z. Coelho, P. Collet, Limit law for the 
close approach of two trajectories in expanding maps
of the circle, Probab. Theory Related Fields 99 (1994), no. 2, 237--250. 

\item{[CC2]} Z. Coelho, P. Collet, Asymptotic limit law for subsystems of 
shifts of finite type, preprint (2000).

\item{[CdF]} Z. Coelho, E. de Faria, Limit laws of entrance times for
homeomorphisms of the circle, Israel Journal of Mathematics 
93 (1996), 93--112.

\item{[CG]} P. Collet, A. Galves, Asymptotic distribution of 
entrance times for expanding maps of the interval, 
Dynamical systems and applications, 139--152, 
World Sci. Ser. Appl. Anal., 4, World Sci. Publishing, 
River Edge, NJ, 1995. 

\item{[DDM]} P. Dartnell, F. Durand, A. Maass, Orbit equivalence
and Kakutani equivalence with Sturmian subshifts, to appear 
in Studia Mathematica.

\item{[D]} F. Durand, Linearly recurrent subshifts have a finite 
number of non-periodic subshift factors,  Ergodic Theory 
and Dynamical Systems 20 (2000), 1061-1078.

\item{[DHS]} F. Durand, B. Host, C. Skau, Substitutive dynamical 
systems, Bratteli diagrams and dimension groups, Ergodic Theory 
and Dynamical Systems 19 (1999), 953--993.

\item{[F]} A. H. Forrest, $K$--groups associated with 
substitution minimal systems,
Israel J. of Math. 98 (1997), 101--139.

\item{[GPS]} T. Giordano, I. Putnam, C. F. Skau, 
Topological orbit equivalence and $C^{*}$--crossed
products, J. reine angew. Math. 469 (1995), 51--111.

\item{[GJ]} R. Gjerde, O. Johansen, 
Bratteli-Vershik models for Cantor minimal systems:
applications to Toeplitz flows, to appear in Ergodic Theory and Dynamical Systems.

\item{[H]} M. Hirata, Poisson law for axiom A diffeomorphisms, Ergodic Theory
and Dynamical Systems 13 (1993), 533--556.

\item{[HSV]} M. Hirata, B. Saussol, S. Vaienti, 
Statistics of return times: a general framework and new
applications, Comm. Math. Phys. 206 (1999), no. 1, 33--55. 

\item{[HW]} G. H. Hardy, E. M. Wright, 
An introduction to the theory of numbers, 4th Edition, 
Oxford (1975).

\item{[HM]}  G. A. Hedlund, M. Morse, Symbolic Dynamics II. Sturmian trajectories, 
American J. of Math. 62, (1940), 1--42.

\item{[HJ]} R. A.  Horn, C. R. Johnson, Matrix Analysis, Cambridge University Press (1985).

\item{[HPS]} R. H. Herman, I. Putnam, C. F. Skau, 
Ordered Bratteli diagrams, dimension groups
and topological dynamics, Internat. J. of Math. 3 (1992), 827--864.

\item{[L]} Y. Lacroix,  Possible limit laws for entrance times 
of an ergodic aperiodic dynamical system, preprint LAMFA-Universit\'e de 
Picardie Jules Verne (2001).

\item{[N]} J. Neveu, Processus Ponctuels, Springer Lecture Notes in Mathematics 598 (1976), 249--445.

\item{[P]} B. Pitskel, Poisson limit law for Markov chains, Ergodic Theory
and Dynamical Systems 11 (1991), 501--513.

\item{[Q]} M. Queff\'elec, Substitution Dynamical Systems, 
Lecture Notes in Mathematics 1294 (1987).

\item{[W]} S. Williams, 
Toeplitz minimal flows which are not uniquely ergodic, 
Z. Wahrsch. Verw. Gebiete 67 (1984), no. 1, 95--107.

\bye